\documentstyle[amssymb,amsfonts, 11 pt]{amsart}

\newtheorem{thm}{Theorem}[section]

\theoremstyle{definition}
\newcommand{\comment}[1]{}

\theoremstyle{remark}

\numberwithin{equation}{section}

\def\be#1{\begin{equation}\label{#1}}

\def \endprf{\hfill  {\vrule height6pt width6pt depth0pt}\medskip}
\def\emph#1{{\it #1}}
\def\textbf#1{{\bf #1}}
\begin{document}       

\title[Navier-Stokes in a critical space]{Ill-posedness of the Navier-Stokes equations in a critical space in 3D}

\author[J. Bourgain]{Jean Bourgain}
\address[J. Bourgain]
{School of Mathematics, Institute for Advanced Study,
1 Einstein Drive, Princeton, New Jersey 08540 US} 
\email{bourgain@@math.ias.edu}

\author[N. Pavlovi\'{c}]{Nata\v{s}a Pavlovi\'{c}}
\address[N. Pavlovi\'{c}]
{Department of Mathematics, University of Texas at Austin,
1 University Station, C1200, Austin, Texas 78712}
\email{natasa@@math.utexas.edu}

\thanks{J.B. was partially supported by NSF grant. 
N.P. was partially supported by NSF grant number DMS 0758247
and by an Alfred P. Sloan Research Fellowship.}

\begin{abstract}
We prove that the Cauchy problem for the three dimensional 
Navier-Stokes equations is ill posed in  $\dot{B}^{-1,\infty}_{\infty}$
in the sense that a ``norm inflation'' happens in finite time. 
More precisely, we show that initial data in the Schwartz class 
$\mathcal{S}$ that are arbitrarily small in $\dot{B}^{-1, \infty}_{\infty}$
can produce solutions arbitrarily large in 
$\dot{B}^{-1, \infty}_{\infty}$ after an arbitrarily short time.
Such a result implies that the solution map itself is discontinuous in 
$\dot{B}^{-1, \infty}_{\infty}$ at the origin. 
\end{abstract}

\maketitle

\section{Introduction}

In this paper we address a long standing open problem concerning 
well-posedness of the three dimensional Navier-Stokes equations in 
the largest critical space  $\dot{B}^{-1,\infty}_{\infty}$ and 
prove that the Cauchy problem for the three dimensional Navier-Stokes equations 
is ill posed in  $\dot{B}^{-1,\infty}_{\infty}$.

The Navier-Stokes equations for the incompressible fluid in ${\mathbb R}^{3}$ are given by
\be{ns} \frac{\partial u}{\partial t} + (u \cdot \nabla) u + \nabla p
= \nu{\Delta}u, 
\end{equation}
\be{nsdiv} \nabla \cdot u = 0,
\end{equation}
and the initial condition
\be{0nsin} u(x,0) = u_{0}(x),
\end{equation}
for the unknown velocity vector $u = u(x,t) \in {\mathbb R}^{3}$ and the pressure $p = p(x,t) \in {\Bbb R}$,
where $x \in {\mathbb R}^{3}$ and $t \in [0,\infty)$. 

We adapt the standard notion of well-posedness. More precisely, 
a Cauchy problem is said to be {\em locally 
well-posed} in $Z$ if for every initial data $u_0(x) \in Z$ 
there exists a time $T = T(\|u_0\|_{Z})>0$ 
such that a solution to the initial value problem exists 
in the time interval $[0,T]$, is unique 
in a certain Banach space of functions $Y \subset C\left([0,T]; Z \right)$
and the solution map from the initial data $u_0$ to the solution $u$ is continuous 
from $Z$ to $C\left( [0,T]; Z \right) $. If $T$ can be taken arbitrarily large 
we say that the Cauchy problem is {\em globally well-posed}.
Also we say that the Cauchy problem is {\em ill-posed} 
if it is not well-posed. Having such a definition of 
ill-posedness it is clear that the problem may be ill-posed due to different reasons ranging 
from a failure of a solution map to be continuous to a more serious type of ill-posedness 
such as a blow-up in finite time. Here we shall establish 
an ill-posedness of the Navier-Stokes initial value problem 
\eqref{ns} - \eqref{0nsin} via proving a finite time blow-up 
for solutions to the Navier-Stokes equations 
in the largest critical space, the Besov space $\dot{B}^{-1, \infty}_{\infty}$. 

In order to understand the role of the space $\dot{B}^{-1, \infty}_{\infty}$ 
in the analysis of the Navier-Stokes
equations we recall the scaling property of the equations first.
It is easy to see that if the pair 
$(u(x,t), p(x,t))$ solves the Navier-Stokes equations \eqref{ns} then 
$(u_{\lambda}(x,t), p_{\lambda}(x,t))$ with
$$u_{\lambda}(x,t) = \lambda \; u(\lambda x, {\lambda}^2 t),$$
$$p_{\lambda}(x,t) = {\lambda}^2 \; p(\lambda x, {\lambda}^2 t)$$
is a solution to the system \eqref{ns} 
with the initial data
$$u_{0 \; \lambda} = \lambda u_0
  (\lambda x) \;\;.$$
The spaces which are invariant under 
such a scaling are called critical spaces for the Navier-Stokes. Examples
of critical spaces for the Navier-Stokes in 3D are: 
\be{embed} 
\dot{H}^{\frac{1}{2}} \hookrightarrow L^3  \hookrightarrow  \dot{B}_{p|p < \infty}^{-1+\frac{3}{p}, \infty}  
\hookrightarrow BMO^{-1} \hookrightarrow \dot{B}^{-1, \infty}_{\infty}.
\end{equation} 
Kato \cite{Ka84} initiated the study of the Navier-Stokes equations in critical spaces 
by proving that the problem \eqref{ns}-\eqref{0nsin} is locally well-posed in $L^3$ 
and globally well-posed if the initial data are small in $L^3({\mathbb R}^3)$.
The study of the Navier-Stokes equations in critical spaces was 
continued by many authors, see, for example, \cite{GM89, T92, C97, P96}. 
In particular, in 2001 Koch and Tataru \cite{KT01} proved 
the global well-posedness of the Navier-Stokes equations evolving 
from small initial data in the space $BMO^{-1}$. 
The space $BMO^{-1}$ has a special role since it is the largest 
critical space among the spaces listed in \eqref{embed} 
where such existence results are available. 

The importance of considering the three dimensional Navier-Stokes equations 
in the Besov space $\dot{B}_{\infty}^{-1,\infty}$ is related to the fact that 
all critical spaces for the 3D Navier-Stokes equations are embedded in the 
same function space, $\dot{B}_{\infty}^{-1,\infty}$. A proof of this
embedding could be found in, for example, \cite{Cannone04}. 
It has been a long standing problem to determine if the Navier-Stokes initial value  
problem is well-posed in the space $\dot{B}^{-1, \infty}_{\infty}$. 
The problem is stated as a conjecture in 
\cite{Cannone04} and \cite{Meyer99}.

An indication that the Navier-Stokes initial value problem might be ill-posed 
in the largest critical space is given in \cite{MS}, where Montgomery-Smith proved 
a finite time blow-up for solutions of a simplified model for the Navier-Stokes 
equations in the space $\dot{B}^{-1, \infty}_{\infty}$. The work \cite{MS} 
suggests that the applications of a fixed point argument that are available up to now
are not likely to produce an existence result for the Navier-Stokes equations 
themselves in the largest critical space, but it does not prove this for the 
actual Navier-Stokes equations.

In this paper we prove that the actual Navier-Stokes system is ill-posed 
in $\dot{B}^{-1,\infty}_{\infty}$ in the sense that there is 
a so called 
``norm inflation'' 
(for similar results in the context of NLS see, e.g. \cite{CCT2003-2}).
Here by a ``norm inflation'' we mean that 
initial data in the Schwartz class $\mathcal{S}$ that are arbitrarily small in $\dot{B}^{-1, \infty}_{\infty}$
can produce solutions arbitrarily large in 
$\dot{B}^{-1, \infty}_{\infty}$ after an arbitrarily short time.
Such a result implies that the solution map itself is discontinuous in 
$\dot{B}^{-1, \infty}_{\infty}$
at
the origin. 
More precisely, our ``norm inflation'' result can be formulated in the following way: 
\begin{thm} \label{main}
For any $\delta > 0$ there exists a solution $(u,p)$ to the Navier-Stokes equations 
\eqref{ns} - \eqref{0nsin} and $0 < t < \delta$ such that $u(0) \in {\cal S}$ 
$$ \|u(0)\|_{\dot{B}^{-1,\infty}_{\infty}} \leq \delta,$$ 
with 
$$\|u(t)\|_{\dot{B}^{-1,\infty}_{\infty}} > \frac{1}{\delta}.$$ 
\end{thm}

We remark that similar programs of establishing ill-posedness have been successfully carried 
out in the context of the nonlinear dispersive equations, see for example work of 
Bourgain \cite{Bo}, Kenig, Ponce, Vega \cite{kpv}, 
Christ-Colliander-Tao \cite{CCT2003-2}, \cite{CCT2003-3}.

The main idea of our approach is to choose initial data $u_0$ in 
$\dot{B}^{-1,\infty}_{\infty} \cap {\mathcal S}$ 
so that when they evolve in time a certain part of the solution  
will become arbitrarily large in finite time. More precisely, 
we write a solution to the Navier-Stokes equations 
\eqref{ns} - \eqref{0nsin} as
$$ 
u = e^{t \Delta} u_0 - u_1 + y, 
$$
where $u_1$ is the first approximation of the
solution to the corresponding linear equation 
and is given by 
$$
u_1(x,t) = \int_{0}^{t} e^{(t-\tau) \Delta} {\mathbb P} 
(e^{\tau \Delta} u_0 \cdot \nabla) e^{\tau \Delta} u_0 \; d\tau,
$$
where ${\mathbb P}$ denotes the projection on divergence free vector fields.
We decompose $u_1$ as $u_1 = u_{1,0} + u_{1,1}$, so that the piece $u_{1,0}$ 
gets arbitrarily large in finite time. 
On the other hand, we obtain a PDE that $y$ solves, thanks to which we control 
$e^{t\Delta}u_0 - u_{1,1} + y$
in the space $X_T$ that was  introduced in \cite{KT01} by Koch and Tataru 
(see Section \ref{prel} for a precise definition of $X_T$).

We note that recently Chemin and Gallagher \cite{CG06} established global existence
of solutions for the Navier-Stokes equations evolving from arbitrary large 
initial data in $\dot{B}^{-1, \infty}_{\infty}$ under the assumption of 
a certain nonlinear smallness on the initial data. Since the initial data that 
we exhibit do not appear to satisfy this nonlinear smallness condition,
our work could be understood as a complement of \cite{CG06}.

After we completed the present paper we learned about the recent 
work of Germain \cite{G08} where he proves 
an instability result for the Navier-Stokes equations in 
$\dot{B}^{-1, q}_{\infty}$, for $q>2$ 
by showing that the map from the initial data to the solution is not 
in the class ${\mathcal{C}}^2$. We remark that \cite{G08} 
does not treat a norm inflation phenomenon.

\subsection*{Organization of the paper} In section \ref{prel} we
introduce the notation that shall be used throughout the paper. 
Also in Section \ref{prel} we recall the result of Koch and Tataru \cite{KT01}. 
In section \ref{pr} we present a proof of Theorem \ref{main}.

\section{Preliminaries} \label{prel}

\subsection{Notation} 

We shall denote by $A \lesssim B$ an estimate of the form $A \leq C B$ with some constant $C$.
Throughout the paper, $i^{\mbox{th}}$ coordinate ($i = 1,2,3$) of 
a vector $x \in {\mathbb R}^3$  will be denoted by $x^i$. 

We recall that the Besov space $\dot{B}^{-1, \infty}_{\infty}$ is equipped 
with the norm 
$$ \|f(\cdot)\|_{\dot{B}^{-1, \infty}_{\infty}} = 
\sup_{t > 0} t^{\frac{1}{2}} \| e^{t\Delta} f(\cdot)\|_{L^{\infty}}.$$

\subsection{The result of Koch and Tataru}

Here we recall the result of Koch and Tataru \cite{KT01} that 
establishes the global well-posedness of the Navier-Stokes equations evolving from small initial data in the space $BMO^{-1}$. 

First, let us recall the definition of the space $BMO^{-1}$ as given in \cite{KT01}:
\be{defbmo-1KT}
\|f(\cdot)\|_{BMO^{-1}} = \sup_{x_0,R} 
\left( \; \frac{1}{|B(x_0,\sqrt{R})|} \int_{0}^{R} \int_{B(x_0,\sqrt{R})}  
|e^{t \Delta} f(y)|^{2} \; dy \; dt \; \right)^{\frac{1}{2}}.
\end{equation} 
In \cite{KT01} Koch and Tataru proved the following existence theorem:
\begin{thm} \label{KTtheorem} 
The Navier-Stokes equations \eqref{ns} - \eqref{0nsin} have a unique global solution in $X$ 
\begin{align*}
\| u(\cdot, \cdot) \|_{ X} & = \sup_{t} t^{\frac{1}{2}} \|u(\cdot, t)\|_{L^{\infty}} \\
& + \sup_{x_0,R} \left( \; \frac{1}{|B(x_0,\sqrt{R})|} \int_{0}^{R} \int_{B(x_0,\sqrt{R})}  
|u(y,t)|^{2} \; dy \; dt \; \right)^{\frac{1}{2}},
\end{align*}
for all initial data $u_0$ with $\nabla \cdot u_0 = 0$ which are small in $BMO^{-1}$.
\end{thm}

Let $T \in (0, \infty]$. We denote by $X_T$ the space equipped with the norm 
\begin{align*}
\| u(\cdot, \cdot) \|_{ X_T} & = \sup_{0 < t < T} t^{\frac{1}{2}} \|u(\cdot, t)\|_{L^{\infty}} \\
& + \sup_{x_0} \sup_{0 < R <T} \left( \; \frac{1}{|B(x_0,\sqrt{R})|} \int_{0}^{R} \int_{B(x_0,\sqrt{R})}  
|u(y,t)|^{2} \; dy \; dt \; \right)^{\frac{1}{2}}.
\end{align*}

Now let ${\mathbb P}$ denote the projection on divergence free vector fields.
As shown in \cite{KT01}, see also \cite{L00}, the bilinear operator 
\be{bilinop} 
{\mathcal B}(u,v) = \int_{0}^{t} e^{(t-\tau) \Delta} {\mathbb P} (u \cdot \nabla) v \; d\tau,
\end{equation} 
maps $X_T \times X_T$ into $X_T$. More precisely, 
\be{bil} 
\|{\mathcal B}(u,v) \|_{X_T} \lesssim \|u\|_{X_T} \|v\|_{X_T}. 
\end{equation}

\setcounter{equation}{0}
\section{Proof of Theorem \ref{main}} \label{pr}

We rewrite the Navier-Stokes equations \eqref{ns} in the following way: 
\be{nsint} 
u = e^{t \Delta} u_0 - u_1 + y, 
\end{equation} 
where 
\be{u_1}
u_1 (x,t) = {\mathcal B}(e^{t \Delta}u_0(x), e^{t \Delta} u_0(x)), 
\end{equation}  
and $y$ satisfies the following equation:
$$ \partial_{t} y - \Delta y + G_1+G_2+G_3 = 0,$$
where
\begin{align}\begin{split}\label{G} 
& G_1 =  {\mathbb P}[\; (e^{t \Delta} u_0 \cdot \nabla) y + (u_1 \cdot \nabla)y + 
(y \cdot \nabla) e^{t \Delta} u_0 + (y \cdot \nabla) u_1 \;] \\
& G_2 =  {\mathbb P}[\; (y \cdot \nabla) y \;]\\
& G_3 = {\mathbb P}[\; (e^{t \Delta} u_0 \cdot \nabla) u_1 + (u_1 \cdot \nabla) e^{t \Delta} u_0 
+ (u_1 \cdot \nabla) u_1 \; ]. 
\end{split} \end{align} 
 
We shall choose initial data $u_0$ in such a way that when they evolve in time,
the part of the solution $u_1$ will become arbitrarily large in $\dot{B}^{-1,\infty}_{\infty}$ at certain time $T$, 
while we will be able to control the behavior of $y$ in the space $X_T$.

\subsection{Choice of initial data} 

Fix small numbers $T > 0$, $\delta > 0$ and a large number $Q>0$ (eventually
$T \rightarrow 0$, $\delta \rightarrow 0$ and $Q \rightarrow \infty$).
Let $\eta \in {\mathbb S}^2$. Let $r = r(Q)$ be a large integer (to be specified). 
We choose the initial data as follows: 
\be{init} 
u_0 = \frac{Q}{\sqrt{r}} \sum_{s=1}^{r} |k_s| [v_s \cos (k_s \cdot x) + v_{s}^{'} \cos (k_{s}^{'} \cdot x)], 
\end{equation} 
where 
\begin{enumerate} 
\item The vectors  $k_s \in {\mathbb{R}}^{3}$ are parallel to a given vector 
$k_0 \in {\mathbb{R}}^{3}$ and $k_{s}^{'} \in {\mathbb{R}}^{3}$ is defined  
by 
\be{et}
k_s - k_{s}^{'} = \eta.
\end{equation} 
Furthermore, we take $|k_0|$ large (depending on $Q$) and $|k_s|$ ($1 \leq s \leq r$)
very lacunary. For example, 
$$|k_s| = 2^s |k_0| \; |k_{s-1}|,\; \; s=2,3,...,r.$$ 

\item $v_s$, $v_{s}^{'} \in {\mathbb S}^2$ such that 
\begin{enumerate} 
\item 
\be{divtobe} 
k_s \cdot v_s = 0 = k_{s}^{'} \cdot v_{s}^{'}.
\end{equation}
Note that \eqref{divtobe} implies that $\mbox{div} \; u_0 = 0$.

\item By \eqref{et} we may ensure that 
$$ v_s \approx v_s^{'} \approx v \in  {\mathbb S}^2.$$
We require that
\be{12eta} 
\eta \cdot v_s = \eta \cdot v_{s}^{'} = \eta \cdot v = \frac{1}{2}.
\end{equation} 
\end{enumerate}
\end{enumerate}  

It is obvious from \eqref{init} that 
$$ \|u_0\|_{\dot{B}^{-1,\infty}_{\infty}} \sim \frac{Q}{\sqrt{r}} < \delta$$
for appropriate $r$.

\subsection{Analysis of ${\mathbf{u_1}}$} 

Now we analyze $u_1$ with a goal to split it into two pieces $u_{1,0}$ and $u_{1,1}$ such that 
the upper bound on $u_{1,0}$ in the Besov space $\dot{B}^{-1, \infty}_{\infty}$ is roughly $Q^2$ on a certain 
time interval.

For the initial data $u_0$ given by \eqref{init}, 
$e^{\tau \Delta} u_0$ can be written as follows
\be{op-u_0} 
e^{\tau \Delta} u_0 =  \frac{Q}{\sqrt{r}} \sum_{s=1}^{r} |k_s| 
\left(v_s \cos(k_s \cdot x) e^{-|k_s|^2 \tau} + v_{s}^{'} \cos (k_{s}^{'} \cdot x) e^{-|k_{s}^{'}|^2\tau} \right). 
\end{equation} 
Hence we can calculate $\left( e^{\tau \Delta}u_0 \cdot \nabla \right) e^{\tau \Delta}u_0$ 
via its coordinates as follows: 
\begin{align} 
\left(
\left( e^{\tau \Delta}u_0 \cdot \nabla \right) e^{\tau \Delta}u_0 
\right)^{i}
& =  \sum_{j} \partial_{j} [ (e^{\tau \Delta} u_0)^i \; (e^{\tau \Delta} u_0)^j] \nonumber \\
& \sim N_1^i + N_2^i + N_3^i, \label{u_1-non1} 
\end{align} 
where 
\begin{align*} 
N_1^i = & \frac{Q^2}{r} \sum_{s=1}^{r} |k_s|^2 e^{-2|k_s|^2 \tau} \sin(\eta \cdot x) 
[ (\eta \cdot v_{s}^{'}) v_s^i + (\eta \cdot v_s) (v_{s}^{'})^i] \\
N_2^i = & \frac{Q^2}{r} \sum_{s=1}^{r} |k_s|^2 e^{-(|k_s|^2+|k_{s}^{'}|^2) \tau} \sin((k_s + k_{s}^{'})\cdot x) \times \\
& \times [ ((k_s + k_{s}^{'})\cdot v_{s}^{'}) v_s^i + ((k_s + k_{s}^{'}) \cdot v_s) (v_{s}^{'})^i] \\
N_3^i = & \frac{Q^2}{r} \sum_{s \not= s'} |k_s| \; |k_{s'}| e^{-(|k_s|^2+|k_{s'}|^2) \tau} \sin((k_s \pm k_{s'})\cdot x) \times \\ 
&  \times [ ((k_s \pm k_{s'})\cdot v_{s'}) v_s^i + ((k_s \pm k_{s'}) \cdot v_s) v_{s'}^i]
+ \mbox{ similar terms }.
\end{align*}

We consider contributions to $u_1$ coming from each of three terms $N_1$, $N_2$, $N_3$.  
Contributions coming from $N_1$ can be estimated by integrating in time 
and using \eqref{12eta} as follows 
\begin{align*} 
& \int_{0}^{t} e^{(t-\tau)\Delta} N_1 \; dt \\
& \sim \frac{Q^2}{r} \sum_{s=1}^{r} |k_s|^2 \; [\int_{0}^{t}e^{-(t-\tau){|\eta|}^2 - 2|k_s|^2 \tau} \; d\tau] \; \sin(\eta \cdot x) \;  
[ (\eta \cdot v_{s}^{'}) v_s + (\eta \cdot v_s) v_{s}^{'}] \\
& \sim Q^2  \; \sin(\eta \cdot x) \; v,  
\end{align*} 
for 
\be{int} 
\frac{1}{|k_1|^2} \ll T \ll 1.
\end{equation} 
Therefore, recalling \eqref{12eta}
\be{bes-N_1}
\|\int_{0}^{t} e^{(t-\tau)\Delta} {\mathbb P}(N_1) \; dt\|_{B_{\infty}^{-1,\infty}} \sim Q^2.
\end{equation} 
Also 
\be{KT-N_1}
\|\int_{0}^{t} e^{(t-\tau)\Delta} {\mathbb P}(N_1) \; dt\|_{X_T} \lesssim \sqrt{T} Q^2.
\end{equation}

Now consider contributions to $u_1$ coming from $N_3$. 
\begin{align} 
& \left| \int_{0}^{t} e^{(t-\tau)\Delta} N_3 \; dt \right| \nonumber \\
\lesssim & \frac{Q^2}{r} \sum_{s=1}^r \sum_{s' < s} |k_s|\;|k_{s'}| 
\; \left| \int_{0}^{t}e^{-(t-\tau)|k_s \pm k_{s'}|^2 - (|k_s|^2+|k_{s'}|^2)\tau} \; d\tau \right| \times \nonumber \\
& \times | \sin((k_s \pm k_{s'})\cdot x) | \; O(|k_s|) \nonumber \\
\sim & \frac{Q^2}{r} \sum_{s=1}^r \sum_{s' < s} |k_s|\;|k_{s'}| 
\left| \frac{e^{-(|k_s|^2 + |k_{s'}|^2)t} - e^{-|k_s \pm k_{s'}|^2t}} 
{|k_s \pm k_{s'}|^2 - (|k_s|^2 + |k_{s'}|^2)} \right| \; O(|k_s|) \nonumber \\
\lesssim & \frac{Q^2}{r} \sum_{s=1}^r \sum_{s' < s} |k_s|\;|k_{s'}| 
e^{-\frac{1}{2}|k_{s'}|^2 t} \; t \; O(|k_s|)\label{bound-e} \\
\lesssim & \frac{Q^2}{r} \sum_{s=1}^r  |k_{s-1}| 
e^{-\frac{1}{l}k_s^2 t}, \label{N_3-cont}
\end{align}
where to obtain \eqref{bound-e} we use the boundedness 
of the function $g(t) = \frac{1 - e^{-\lambda t}}{\lambda t}$, with $\lambda > 0$, 
while to obtain \eqref{N_3-cont} we use the boundedness of the function 
$h(t) = \mu t e^{-\mu t}$, with $\mu > 0$ 
and we replace $e^{-k_{s'}^2}$ by $e^{-\frac{1}{l} k_s^2}$ for some $l$.
We also use the lacunarity of the sequence $|k_s|$. 

Thus \eqref{N_3-cont} implies that 
\begin{align}
& \|\int_{0}^{t} e^{(t-\tau)\Delta} {\mathbb P}(N_3) \; dt\|_{X_T} \nonumber \\
& \lesssim \frac{Q^2}{r} \sum_{s=1}^r \frac{|k_{s-1}|}{|k_s|} 
+ \frac{Q^2}{r} \sup_{t<T} 
\left\{ \int_0^t 
\left[ \sum_{s=1}^r |k_{s-1}| e^{-\frac{1}{l}|k_s|^2 \tau} \right]^2 \; d\tau 
\right\}^{ \frac{1}{2} } \nonumber \\ 
& \lesssim \frac{Q^2}{r} \sum_{s=1}^r \frac{|k_{s-1}|}{|k_s|} \nonumber \\
& < \frac{Q^2}{r}, \label{KT-N_3}
\end{align} 
again by lacunarity of $|k_s|$.

Next we estimate the contribution coming from $N_2$. 
Clearly, recalling \eqref{et} 
$$ \int_{0}^{t} e^{(t-\tau)\Delta} N_2 \; dt 
\sim \frac{Q^2}{r} 
\left\{ \sum_{s=1}^r O(|k_s| e^{-|k_s|^2 t}) \sin(k_s + k_{s}^{'}) \cdot x
\right\}.$$ 
Therefore 
\begin{align}
& \|\int_{0}^{t} e^{(t-\tau)\Delta} {\mathbb P}(N_2) \; dt\|_{X_T} \nonumber \\
& \lesssim \frac{Q^2}{r} \sup_{t > 0}
\left| \sum_{s = 1}^{r} t^{\frac{1}{2}} |k_s| e^{-|k_s|^2 t} \right| \nonumber \\
& + \frac{Q^2}{r} \sup_{R>0} 
\left\{ 
\int_{0}^{R} 
\left[ \sum_{|k_s| > \frac{1}{\sqrt{R}}} |k_s|^2 e^{-|k_s|^2 t} \right] \; dt 
+ \int_{0}^{R} \left( \sum_{|k_s| \leq \frac{1}{\sqrt{R}}} |k_s| \right)^2 \; dt  
\right\}^{\frac{1}{2}} \nonumber \\
& \lesssim 
\frac{Q^2}{r} + \frac{Q^2}{r} (r + 1)^{\frac{1}{2}} \nonumber \\ 
& \lesssim \frac{Q^2}{\sqrt{r}}, \label{KT-N_2}
\end{align} 
using the fact that $ \sqrt{t} \sum_s |k_s|^2 e^{-|k_s|^2 t} \lesssim 1$ and
making the appropriate splitting to bound the second term in $ \| \cdot \|_{X^T}$.

Hence we can decompose $u_1$ as follows
$$u_1 = u_{1,0} + u_{1,1},$$
where
\begin{align} \begin{split} \label{u1est} 
& \|u_{1,0}\|_{B_{\infty}^{-1,\infty}} \sim Q^2 \mbox{ for } \frac{1}{|k_1|^2} \ll t \ll 1,\\
& \|u_{1,0}\|_{X_T} \lesssim \sqrt{T} Q^2,\\
& \|u_{1,1}\|_{X_T} \lesssim \frac{Q^2}{\sqrt{r}}. 
\end{split} \end{align}

\subsection{Analysis of ${\mathbf{y}}$}
Now we analyze the remaining part of the solution, which we denoted by $y$. 
The main idea is to control $y$ using the space of Koch and Tataru $X_T$. 

Consider time-intervals 
$$0 < T_1 < T_2 < ... < T_{\beta}, \; \; \; \; \beta = Q^3$$
with 
\begin{align} 
& T_{\alpha}^{-1} = |k_{r_{\alpha}}|^2 \label{times-r}\\
& r_{\alpha} = r - \alpha Q^{-3} r, \; \alpha = 1,2, ....\label{r} 
\end{align} 
In particular, $r_{\beta} = 0$ and $T_{\beta}^{-1} = |k_0|^2$.  

For $t \geq T_{\alpha}$ the equation for $y$ can be written in the integral form as 
\be{yDuh}
y(t) = e^{(t-T_{\alpha})\Delta} y(T_{\alpha}) - \int_{T_{\alpha}}^{t} e^{(t-\tau)\Delta} 
[G_1 + G_2 + G_3](\tau) \; d\tau, 
\end{equation}
where $G_i$, $i = 1,2,3$ are given by \eqref{G}. 

Also 
$$y(T_{\alpha}) = \int_{0}^{T_{\alpha}} e^{(T_{\alpha}-\tau)\Delta} 
[G_1 + G_2 + G_3](\tau) \; d\tau.$$ 
Therefore 
\begin{align} 
e^{(t-T_{\alpha})\Delta} y(T_{\alpha}) 
& = \int_{0}^{T_{\alpha}} e^{(t-\tau)\Delta} [G_1 + G_2 + G_3](\tau) \; d\tau \nonumber \\
& = \int_{0}^{t} e^{(t-\tau)\Delta} [G_1 + G_2 + G_3](\tau) \; \chi_{[0,T_{\alpha}]} \; d\tau, \label{op-yTalp}
\end{align}
where  $\chi_{[0,T_{\alpha}]}$ is a characteristic function of the interval  $[0,T_{\alpha}]$. 

Now we substitute \eqref{op-yTalp} in \eqref{yDuh} to obtain 
\be{y-Tal1} 
\|y\|_{ X_{T_{\alpha +1}} } \leq I + II, 
\end{equation} 
where  
\be{0Tal} 
I = \|\int_{0}^{t} e^{(t-\tau)\Delta} [G_1 + G_2 + G_3](\tau) \; \chi_{[0, T_{\alpha}]}(\tau) \; d\tau 
\|_{ X_{T_{\alpha +1}} } 
\end{equation}
and 
\be{Tal} 
II = \|\int_{0}^{t} e^{(t-\tau)\Delta} [G_1 + G_2 + G_3](\tau) \; \chi_{[T_{\alpha}, T_{\alpha +1}]}(\tau) \; d\tau 
\|_{ X_{T_{\alpha +1}} }.
\end{equation}

Next we use the bilinear estimate \eqref{bil} 
on the terms in $G_1$, $G_2$ and $G_3$ to obtain 
an upper bound on $I$ and $II$ respectively. 
Before we obtain an upper bound on $I$, we estimate
$\| e^{t\Delta} u_0 \|_{ X_{T_{\alpha}} }$. 
From \eqref{op-u_0} we have 
$$
e^{t \Delta} u_0 \approx
\frac{Q}{\sqrt{r}} \sum_{s \leq r} |k_s| 
v_s \cos(k_s \cdot x) e^{-|k_s|^2 t}.
$$  
We estimate $\|e^{t \Delta} u_0\|_{X_{T_{\alpha}}}$ as follows 
\begin{align}
\| e^{t \Delta} u_0 \|_{X_{T_{\alpha}}}  
\leq  & \frac{Q}{\sqrt{r}} \sup_{t < T_{\alpha}} \sqrt{t} \sum_{ s \leq r} |k_s| e^{-k_s^2 t}  \label{KTal-u0-first} \\  
& + \frac{Q}{\sqrt{r}} \sup_{x_0,\; \; 0 < t < T_{\alpha}} \left( t^{-3/2} \int_{0}^{t} \int_{|x-x_0| < \sqrt{t}} 
\left| \sum_{ s \leq r} 
|k_s| v_s \cos(k_s \cdot x) e^{-k_s^2 \tau} \right|^2 
\; dx \; d\tau \right)^{1/2} \nonumber \\
\lesssim  & \frac{Q}{\sqrt{r}} + 
\frac{Q}{\sqrt{r}} (r+1)^{\frac{1}{2}} \nonumber \\
\lesssim & Q,  
\end{align} 
similarly to \eqref{KT-N_2}. 
Hence 
\be{KTop-u0-al}
\| e^{t\Delta} u_0 \|_{ X_{T_{\alpha}} } \lesssim
Q.
\end{equation}
Now we are ready to estimate $I$ using \eqref{G} and the bilinear 
estimate \eqref{bil}:
\begin{align} \begin{split} \label{0Talest} 
I \lesssim & 
\left( \| e^{t\Delta} u_0 \|_{ X_{T_{\alpha}} } 
+ \| u_1 \|_{ X_{T_{\alpha}} } 
+ \|y\|_{ X_{T_{\alpha}} } \right) 
\|y\|_{ X_{T_{\alpha}} } \\
& + \left( \| e^{t\Delta} u_0 \|_{ X_{T_{\alpha}} } 
+ \| u_1 \|_{ X_{T_{\alpha}} } \right)
\| u_1 \|_{ X_{T_{\alpha}} } \\
\leq  & \left(Q + Q^2 T_{\alpha}^{1/2} + \frac{Q^2}{\sqrt{r}} + \|y\|_{ X_{T_{\alpha}} } \right) \|y\|_{ X_{T_{\alpha}} } \\
& + \left( Q + Q^2 T_{\alpha}^{1/2} + \frac{Q^2}{\sqrt{r}} \right)
\left( Q^2 T_{\alpha}^{1/2} + \frac{Q^2}{\sqrt{r}} \right), 
\end{split} \end{align}
where to obtain \eqref{0Talest} we used \eqref{KTop-u0-al} and \eqref{u1est}.

In order to obtain an upper bound on $II$, first, we estimate 
$\| (e^{t\Delta} u_0) \; \chi_{[T_{\alpha}, T_{\alpha +1}]} (t) \|_{ X_{T_{\alpha +1}} }$. 
More precisely, from \eqref{op-u_0} we have 
\be{linal1}
( e^{t \Delta} u_0) \; \chi_{[T_\alpha, T_{\alpha +1}]}(t) \approx L_1 + L_2, 
\end{equation}
where 
$$ L_1 =  \frac{Q}{\sqrt{r}} \sum_{s < r_{\alpha +1}} |k_s| v_s \cos(k_s \cdot x) 
\chi_{[T_{\alpha}, T_{\alpha +1}]}(t)$$ 
and 
$$L_2 =  \frac{Q}{\sqrt{r}} \sum_{s=r_{\alpha +1} }^{r_{\alpha}} |k_s| 
v_s \cos(k_s \cdot x) e^{-|k_s|^2 t} \chi_{[T_{\alpha}, T_{\alpha +1}]}(t).$$  
We estimate $L_1$ keeping in mind that, thanks to \eqref{times-r},
$T_{\alpha + 1} = |k_{r_{\alpha+1}}|^{-2}$:   
\begin{align}
\| L_1 \|_{ X_{T_{\alpha +1}} } 
\leq  & \frac{Q}{\sqrt{r}} T_{\alpha + 1}^{1/2} |k_{ r_{\alpha + 1} - 1 }| \nonumber \\
& + \frac{Q}{\sqrt{r}} \sup_{x_0, t} \left( t^{-3/2} \int_{0}^{t} \int_{|x-x_0| < \sqrt{t}} 
\left| \sum_{ s < r_{\alpha +1} } |k_s| \chi_{[T_{\alpha}, T_{\alpha+1}]}(\tau) \right|^2 \; dx \; d\tau \right)^{1/2} \nonumber \\
\leq & \frac{Q}{\sqrt{r}} \frac{  |k_{ r_{\alpha + 1} - 1 }| }{ |k_{r_{\alpha + 1}}| } + 
 \frac{Q}{\sqrt{r}} \left( T_{\alpha + 1} |k_{ r_{\alpha + 1} - 1 }|^2 \right)^{1/2} \nonumber \\
< & \frac{Q}{\sqrt{r}}. \label{L1} 
\end{align} 
We estimate $L_2$ as follows 
\begin{align}
\| L_2 \|_{X_{T_{\alpha +1}} } 
\leq  & \frac{Q}{\sqrt{r}} \sup_{t} \sqrt{t} \sum_{ s = r_{\alpha+1} }^{r_{\alpha}} |k_s| e^{-k_s^2 t} \nonumber \\  
& + \frac{Q}{\sqrt{r}} \sup_{x_0, t} \left( t^{-3/2} \int_{0}^{t} \int_{|x-x_0| < \sqrt{t}} 
\left| \sum_{ s = r_{\alpha +1} }^{r_{\alpha}} 
|k_s| v_s \cos(k_s \cdot x) e^{-k_s^2 \tau} \chi_{[T_{\alpha}, T_{\alpha+1}]}(\tau) \right|^2 
\; dx \; d\tau \right)^{1/2} \nonumber \\
\leq & \frac{Q}{\sqrt{r}} +  \frac{Q}{\sqrt{r}} (r_{\alpha} - r_{\alpha +1})^{1/2} \nonumber \\
= & \frac{Q}{\sqrt{r}} +  \frac{Q}{\sqrt{r}} (Q^{-3} r)^{1/2} \label{use-r} \\
\leq & Q^{-1/2}, \label{L2}  
\end{align} 
where to obtain \eqref{use-r} we used \eqref{r}.  
Hence we combine \eqref{linal1}, \eqref{L1} and \eqref{L2} to conclude 
\be{KTop-u0}
\| (e^{t\Delta} u_0) \; \chi_{[T_{\alpha}, T_{\alpha +1}]} (t) \|_{ X_{T_{\alpha +1}} } \lesssim
Q^{-1/2}.
\end{equation}

Also we recall that \eqref{u1est} implies 
\be{KTu1}
\| u_1 \; \chi_{[T_{\alpha}, T_{\alpha +1}]} (t) \|_{ X_{T_{\alpha +1}} } \lesssim
Q^2 T_{\alpha + 1}^{1/2} + \frac{Q^2}{\sqrt{r}}.
\end{equation} 

Now we are ready to find an upper bound on $II$ by employing the bilinear estimate \eqref{bil}: 
 \begin{align} \begin{split} \label{Talest} 
II 
\lesssim & \left( \| (e^{t\Delta} u_0) \; \chi_{[T_{\alpha}, T_{\alpha +1}]} (t) \|_{ X_{T_{\alpha +1}} } 
+ \| u_1 \; \chi_{[T_{\alpha}, T_{\alpha +1}]} (t) \|_{ X_{T_{\alpha +1}} } 
+ \|y\|_{ X_{T_{\alpha +1}} } \right) 
\|y\|_{ X_{T_{\alpha +1}} } \\
& + \left(\| (e^{t\Delta} u_0) \; \chi_{[T_{\alpha}, T_{\alpha +1}]} (t) \|_{ X_{T_{\alpha +1}} } 
+ \| u_1 \; \chi_{[T_{\alpha}, T_{\alpha +1}]} (t) \|_{ X_{T_{\alpha +1}} } \right)
\| u_1 \; \chi_{[T_{\alpha}, T_{\alpha +1}]} (t) \|_{ X_{T_{\alpha +1}} } \\
\leq &   
\left( Q^{-1/2} + Q^2 T_{\alpha + 1}^{1/2} + \frac{Q^2}{\sqrt{r}} + \|y\|_{ X_{T_{\alpha +1}} } \right)
\|y\|_{ X_{T_{\alpha +1}} } \\
& + \left( Q^{-1/2} + Q^2 T_{\alpha + 1}^{1/2} + \frac{Q^2}{\sqrt{r}} \right) 
\left( Q^2 T_{\alpha + 1}^{1/2} + \frac{Q^2}{\sqrt{r}} \right), 
\end{split} \end{align} 
where to obtain \eqref{Talest} we used \eqref{KTop-u0} and  \eqref{KTu1}.

Having in mind that $T_{\alpha} < T_{\alpha+1} < T$ and that $T$ 
will be chosen to satisfy \eqref{choice},
we combine 
\eqref{y-Tal1}, \eqref{0Talest} and \eqref{Talest} to obtain
$$ \| y \|_{ X_{T_{\alpha +1}} } \lesssim 
Q^{-1/2} \|y\|_{ X_{T_{\alpha +1}} }
+ \|y\|_{ X_{T_{\alpha +1}} }^{2}
+ Q^3 (\frac{1}{\sqrt{r}} + T_{\alpha +1}^{1/2}) 
+ Q \|y\|_{ X_{T_{\alpha}} }.$$ 
Thus 
\be{yal1}
\| y \|_{ X_{T_{\alpha +1}} } \lesssim 
Q^3 (\frac{1}{\sqrt{r}} + T_{\beta}^{1/2}) + Q \|y\|_{ X_{T_{\alpha}} }.
\end{equation} 
Iterating \eqref{yal1} gives 
\be{ybet}
\| y \|_{ X_{T_{\beta}} } \lesssim 
Q^{\beta+3} (\frac{1}{r} + T_{\beta})^{1/2}. 
\end{equation} 

Now we take $T>T_{\beta}$ and write \eqref{yDuh} and \eqref{op-yTalp} with $\alpha = \beta$. 
Thus
\be{y-T}  
\|y\|_{ X_{T} } \leq I_{\beta} + II_{\beta} 
\end{equation} 
where 
\be{0Tbe} 
I_{\beta} = \|\int_{0}^{t} e^{(t-\tau)\Delta} [G_1 + G_2 + G_3](\tau) \; \chi_{[0, T_{\beta}]}(\tau) \; d\tau 
\|_{ X_{T} } 
\end{equation} 
and 
\be{Tbe} 
II_{\beta} = \|\int_{0}^{t} e^{(t-\tau)\Delta} [G_1 + G_2 + G_3](\tau) \; \chi_{[T_{\beta}, T]}(\tau) \; d\tau 
\|_{ X_{T} }.
\end{equation} 

We obtain an upper bound on $I_{\beta}$ by using \eqref{G} and 
the bilinear estimate \eqref{bil}:
\begin{align}  
I_{\beta} 
\lesssim &
\left( \| (e^{t\Delta} u_0) \|_{ X_{T_{\beta}} } 
+ \| u_1 \|_{ X_{T_{\beta}} } 
+ \|y\|_{ X_{T_{\beta}} } \right) 
\|y\|_{ X_{T_{\beta}} } \nonumber \\
& + \left( \| (e^{t\Delta} u_0) \|_{ X_{T_{\beta}} } 
+ \| u_1 \|_{ X_{T_{\beta}} } \right)
\| u_1 \|_{ X_{T_{\beta}} } \nonumber \\
\leq & 
\left(Q + Q^2 T_{\beta}^{1/2} + \frac{Q^2}{\sqrt{r}} +  
Q^{Q^3} (\frac{1}{r} + T_{\beta})^{1/2}\right) 
Q^{Q^3} (\frac{1}{r} + T_{\beta})^{1/2} \label{0Tbe1}  \\
& + \left( Q + Q^2 T_{\beta}^{1/2} + \frac{Q^2}{\sqrt{r}} \right)
\left( Q^2 T_{\beta}^{1/2} + \frac{Q^2}{\sqrt{r}} \right). \label{0Tbe2} 
\end{align}
We rely here on \eqref{KTop-u0-al}, \eqref{u1est} and \eqref{ybet}.

Recalling that $T_{\beta} = |k_0|^{-2}$ and choosing $r$ and $|k_0|$ large
enough, it follows from \eqref{0Tbe1} and \eqref{0Tbe2} that 
\begin{equation} \label{0Tbeest} 
I_{\beta} \lesssim r^{-1/3} + |k_0|^{-1/2}.
\end{equation}

Also 
\begin{align} 
II_{\beta} \lesssim &
\left( \| (e^{t\Delta} u_0) \; \chi_{[T_{\beta}, T]} (\tau) \|_{ X_{T} } 
+ \| u_1 \|_{ X_{T} } 
+ \|y\|_{ X_{T} } \right) 
\|y\|_{ X_{T} } \nonumber \\
& + \left(\| (e^{t\Delta} u_0) \; \chi_{[T_{\beta}, T]} (\tau) \|_{ X_{T} } 
+ \| u_1 \|_{ X_{T} } \right)
\| u_1 \|_{ X_{T} } \nonumber \\
\lesssim &
\left( |k_1| e^{-\frac{|k_1|^2}{|k_0|^2}} 
+ Q^2 T^{1/2} + \frac{Q^2}{\sqrt{r}} + \|y\|_{ X_{T} } \right)
\|y\|_{ X_{T} }
+ \left( |k_1| e^{-\frac{|k_1|^2}{|k_0|^2}}
+ Q^2 T^{1/2} + \frac{Q^2}{\sqrt{r}} \right)^{2} \label{Tbeest}
\end{align} 
where to obtain \eqref{Tbeest} we used \eqref{op-u_0} and \eqref{u1est}.

Let us also assume that 
\be{choice} 
T < Q^{-8}.
\end{equation}  
Since $|k_1| > |k_0|^2$, \eqref{Tbeest} implies 
\be{Tbeest1} 
II_{\beta} < \left(o(1) + \|y\|_{X_T}\right) \|y\|_{X_T} + 2 Q^4 T. 
\end{equation} 
Therefore, from \eqref{0Tbeest} and  \eqref{Tbeest1}  
\be{y-XT} 
\|y\|_{ X_{T} } < 3 Q^4 T
\end{equation} 
implying 
\be{y-inf} 
\|y\|_{L^{\infty}} \leq T^{-\frac{1}{2}} \|y\|_{X_T} <  3 Q^4 T^{\frac{1}{2}}. 
\end{equation}  

Now we combine \eqref{nsint}, \eqref{u1est} and  \eqref{y-inf} to conclude that 
$$ 
\|u(T) - e^{T\Delta}u_0\|_{\dot{B}^{-1, \infty}_{\infty} } 
\geq Q^2 - \|u_{1,1}\|_{L^{\infty}} - \|y\|_{L^{\infty}}  
> Q^2 (1 - \frac{1}{\sqrt{rT}}  - 3Q^2 T^{\frac{1}{2}})
$$
and 
\be{uTBes}  
\|u(T)\|_{\dot{B}^{-1, \infty}_{\infty} } > \frac{1}{2} Q^2.
\end{equation} 

Consequently we proved that for all $\delta > 0$ 
\be{large} 
\sup_{ \|u(0)\|_{\dot{B}^{-1,\infty}_{\infty} } \leq \delta } \; \; 
\sup_{0 < t < \delta} 
\|u(t)\|_{\dot{B}^{-1,\infty}_{\infty} } 
= \infty.
\end{equation}

\end{document}